\documentclass{gtart_h}
\def\ifplaintex{\expandafter\ifx\csname documentclass\endcsname\relax}

\def\gtp{{\mathsurround=0pt\it $\cal G\mskip-2mu$eometry \&\ 
$\cal T\!\!$opology $\cal P\!$ublications}}  % GT publications

\def\recd{{\small Received:\qua\receiveddate\ifx\reviseddate\relax
\else\qquad Revised:\qua\reviseddate\fi\par}} 

\def\lognumber#1{\def\thelognumber{#1}}
\def\volumenumber#1{\def\thevolumenumber{#1}}
\def\volumeyear#1{\def\thevolumeyear{#1}}
\def\papernumber#1{\def\thepapernumber{#1}}
\def\pagenumbers#1#2{\def\startpage{#1}\def\finishpage{#2}}
\def\published#1{\def\publishdate{#1}}

\def\received#1{\def\receiveddate{#1}}

\def\accepted#1{\def\accepteddate{#1}}

\long\def\asciiabstract#1{\long\def\theasciiabstract{#1}}

\let\\\par\let\thelognumber\relax\let\thevolumenumber\relax
\let\thepapernumber\relax\let\thevolumeyear\relax\let\startpage\relax
\let\finishpage\relax\let\publishdate\relax\let\receiveddate\relax
\let\reviseddate\relax\let\accepteddate\relax\let\theasciititle\relax
\let\theasciiauthors\relax
\let\theasciiabstract\relax

\let\theasciiemail\relax

\ifplaintex
\font\logobig=cmssbx10 scaled 3836
\font\logomed=cmssbx10 scaled 2557
\else
\font\logobig=cmssbx10 scaled 4200
\font\logomed=cmssbx10 scaled 2800
\fi

\long\def\makeagttitle{   %%% start of definition of \makeagttitle
\count0=\startpage
\agt\hfill      %   Journal title (top left) 
\hbox to 45truept{\vbox to 0pt{\vglue -13truept{\logomed A\kern -.37em{\logobig 
T}\kern -.38em G}\vss}\hss}
\break
{\small Volume \thevolumenumber\ (\thevolumeyear)
\startpage--\finishpage\nl
Published: \publishdate}

\vglue .25truein

{\parskip=0pt\leftskip 0pt plus
1fil\def\\{\par\smallskip}{\Large\bf\thetitle}\par\medskip} \vglue
0.05truein

{\parskip=0pt\leftskip 0pt plus 1fil\def\\{\par}{\sc\theauthors}
\par\medskip}%
 
\vglue 0.03truein 

{\small\leftskip 25truept\rightskip 25truept{\bf Abstract}\stdspace\theabstract

{\bf AMS Classification}\stdspace\theprimaryclass
\ifx\thesecondaryclass\relax\else; \thesecondaryclass\fi\par
{\bf Keywords}\stdspace \thekeywords\par}\vglue 7truept

}   %%%% end of definition of \makeagttitle

\ifplaintex
\font\phead=cmsl9 scaled 950
\font\pnum=cmbx10 scaled 913
\font\pfoot=cmsl9 scaled 950
\headline{\vbox to 0pt{\vskip -4.5mm\line{\small\phead\ifnum
\count0=\startpage ISSN 1472-2739 (on-line) 1472-2747 (printed)
\hfill {\pnum\folio}\else\ifodd\count0\def\\{ }% 
\ifx\theshorttitle\relax\thetitle\else\theshorttitle\fi\hfill{\pnum\folio}
\else\def\\{ and }{\pnum\folio}\hfill\ifx\theshortauthors\relax\theauthors
\else\theshortauthors\fi\fi\fi}\vss}}
\footline{\vbox to 0pt{\vglue 0mm\line{\small\pfoot\ifnum\count0=\startpage
\copyright\ \gtp\hfill\else
\agt, Volume \thevolumenumber\ (\thevolumeyear)\hfill\fi}\vss}}
\else
\font=cmsl9 scaled 1050
\font\lnum=cmbx10 
\font=cmsl9 scaled 1050
\makeatletter
\def\@oddhead{{\small\lhead\ifnum\count0=\startpage ISSN 1472-2739 
(on-line) 1472-2747 (printed)\hfill {\lnum\number\count0}\else\ifodd\count0
\def\\{ }\ifx\theshorttitle\relax \thetitle \else\theshorttitle\fi\hfill
{\lnum\number\count0}\else\def\\{ and }{\lnum\number\count0}
\hfill\ifx\theshortauthors\relax 
\theauthors\else\theshortauthors\fi\fi\fi}}\def\@evenhead{\@oddhead}
\def\@oddfoot{\small\lfoot\ifnum\count0=\startpage\copyright\ \gtp\hfill\else
\agt, Volume \thevolumenumber\ (\thevolumeyear)\hfill\fi}
\def\@evenfoot{\@oddfoot}
\makeatother
\fi
\let\maketitlepage\makeagttitle
\let\makeshorttitle\maketitlepage
\let\maketitle\maketitlepage

\newwrite\gtoutfile
\long\gdef\makeheadfile{  %%% start of definition of \makeheadfile
{\def\\{, }\def\s{ }
\immediate\openout\gtoutfile head.xxx
\immediate\write\gtoutfile{Proxy-for: \ifx\theasciiauthors\relax
\theauthors\else\theasciiauthors\fi\s<\ifx\theasciiemail\relax\theemail\else\theasciiemail\fi>}
\immediate\write\gtoutfile{\noexpand\\}
\immediate\write\gtoutfile{Authors: \ifx\theasciiauthors\relax
\theauthors\else\theasciiauthors\fi}
{\def\\{ }\immediate\write\gtoutfile{Title: \ifx\theasciititle\relax
\thetitle\else\theasciititle\fi}}
\immediate\write\gtoutfile{Subj-class: GT or SG, GR etc}
\immediate\write\gtoutfile{MSC-class: \theprimaryclass\ifx\thesecondaryclass\relax\else, \thesecondaryclass\fi}
\immediate\write\gtoutfile{Journal-ref: Algebr. Geom. Topol. \thevolumenumber\s
(\thevolumeyear) \startpage-\finishpage}
\immediate\write\gtoutfile{Comments: Published by Algebraic and
Geometric Topology at}
\immediate\write\gtoutfile{\s\s\s  http://www.maths.warwick.ac.uk/agt/AGTVol\thevolumenumber/agt-\thevolumenumber-\thepapernumber.abs.html}
\immediate\write\gtoutfile{\noexpand\\}
\immediate\write\gtoutfile{}
\ifx\theasciiabstract\relax
\immediate\write\gtoutfile{\theabstract}\else
\immediate\write\gtoutfile{\theasciiabstract}\fi
\immediate\write\gtoutfile{}
\immediate\write\gtoutfile{\noexpand\\}
\immediate\write\gtoutfile{}
\immediate\closeout\gtoutfile}}  %%% end of definition of \makeheadfile

\def\maketitlepage{\makeagttitle\makeheadfile}
\let\makeshorttitle\maketitlepage
\let\maketitle\maketitlepage

\lognumber{8}
\volumenumber{5}
\volumeyear{2005}
\papernumber{8}
\pagenumbers{129}{134}
\received{21 September 2004} 
\accepted{25 February 2005}
\published{10 March 2005}

\usepackage{amsmath,amssymb,graphicx,psfrag}

\def\psfraga <#1,#2> #3#4{%
\psfrag {#3}{\smash{\rlap{\kern #1 \raise #2\hbox{#4}}}}}

\newtheorem{thm}{Theorem}[section] 
\newtheorem{lem}[thm]{Lemma}          
\newtheorem{prop}[thm]{Proposition}

\theoremstyle{definition}
\newtheorem{defn}[thm]{Definition}    
\newtheorem*{quest}{Question}

\newcommand{\NS}{MR1329042}%MR96c:20066}
\newcommand{\Efftpft}{MR1950887}%MR2003i:20067}
\newcommand{\Enonhopf}{MR2022477}
\newcommand{\Epatterns}{MR1974064}%MR2004b:20055}
\newcommand{\Elooppaper}{MR2026832}

\newcommand{\RebbechiThesis}{RebbechiThesis}

\newcommand{\fftp}{falsification by fellow traveler property}
\newcommand{\gset}{generating set}
\newcommand{\cg}{Cayley graph}

\begin{document}

\title{Regular geodesic languages
  and the\\falsification by fellow traveler property}                    
\shorttitle{Regular languages and the fellow traveler property}
\authors{Murray Elder}                  
\address{School of Mathematics and
Statistics, University of St Andrews\\North Haugh, St Andrews, Fife,
KY16 9SS, Scotland}
\email{murray@mcs.st-and.ac.uk}
\urladdr{http://www-groups.mcs.st-andrews.ac.uk/~murray/}

\begin{abstract}   
We furnish an example of a finite generating set for a group that does
 not enjoy the \fftp, while the full language of geodesics is regular.
\end{abstract}

\asciiabstract{%
We furnish an example of a finite generating set for a group that does
not enjoy the falsification by fellow traveler property, while the
full language of geodesics is regular.}


\maketitle

\section{Introduction}

In this short note we answer the following question of Neumann and
Shapiro from \cite{\NS}:
\begin{quest}
 Can one find a monoid generating set $A$ of a group $G$ so
that the language of geodesics is regular but $A$ does not have the
\fftp?
\end{quest}
The converse to this statement is Proposition 4.1 of their paper,
which states that if $A$ has the \fftp\ then the full language of
geodesics on $A$ is regular, and this fact is the reason for the
property's existence.  Several authors have used the \fftp\ as a route
to finding other (geometric) properties  of groups; 
Rebbechi uses a version of the property 
to prove that relatively hyperbolic groups are biautomatic \cite{\RebbechiThesis},
and the author exploits the
property to prove that a certain class of groups is almost convex
\cite{\Enonhopf}.  The author discusses various attributes and
extensions of the property in \cite{\Efftpft, \Elooppaper,
\Epatterns}.

In this article we answer the question via an example first given by
Cannon to demonstrate that a group may have a regular language of
geodesics with respect to one generating set but not another. Neumann
and Shapiro include it in \cite{\NS} to prove that the \fftp\ is
generating set dependent.

The author wishes to thank Walter Neumann and Jon McCammond for 
ideas and help with this paper.

\section{Definitions}

\begin{defn}[\rm(Finite state automaton; regular language)]
Let $A$ be a finite set of letters, and let $A^*$ be the set of all
finite strings, including the empty string, that can be formed from
the letters of $A$.  A {\em finite state automaton} is a quintuple
$(S,A,\tau,Y,s_0)$, where $S$ is a finite set of {\em states}, $\tau$
is a map $\tau:S\times A\rightarrow S$, $Y\subseteq S$ are the {\em
accept states}, and $s_0\in S$ is the {\em start state}. A finite
string $w\in A^*$ is {\em accepted} by the finite state automaton if
starting in the state $s_0$ and changing states according to the
letters of $w$ and the map $\tau$, the final state is in $Y$. The set
of all finite strings that are accepted by a finite state automaton is
called the {\em language} of the automaton.  A language $L\subseteq
A^*$ is {\em regular} if it is the language of a finite state
automaton.
\end{defn}

Suppose $G$ is a group with finite generating set $A$. A word in $A^*$ represents a path
 in the \cg\ based at any vertex. Define $d(a,b)$ to be the distance between two points $a$ 
 and $b$ in the \cg\ with respect to the path metric.
 Paths can be parameterized  by non-negative $t\in\mathbb R$
  by defining $w(t)$ as the point at distance $t$ along the path $w$ if $t$ is between $0$ and 
  the length of $w$, and the endpoint of $w$ otherwise.

\begin{defn}[\rm(The (asynchronous) fellow traveler property)]
Paths $u$ and $v$ are said to {\em $k$-fellow travel} if $d(u(t),v(t))\leq k$ for all $t\geq 0$.
They  {\em asynchronously $k$-fellow travel} if there is a non-decreasing proper continuous 
function $\phi:[0,\infty)\rightarrow [0,\infty)$ such that $d(u(t),v(\phi(t)))\leq k$. 
A language $L\subseteq A^*$ enjoys the {\em (asynchronous) fellow traveler property} 
if there is a constant $k$ such that for each $u,v\in L$ that start at the identity and end at distance 
$0$ or $1$ apart in the \cg,   $u$ and $v$ (asynchronously) $k$-fellow travel.
\end{defn}

\begin{defn}[\rm(The (asynchronous) \fftp)]
A finite generating set $A$ for a group $G$ has the {\em  (asynchronous) \fftp} if
there is a constant $k$ such that every non-geodesic word in the
Cayley graph of $G$ with respect to $A$ is (asynchronously) $k$-fellow traveled by a
shorter word.
\end{defn}
The property arises naturally in the context of geodesic regular
languages, and the proof of Proposition 4.1 in \cite {\NS} uses the
property to build an appropriate finite state automaton. The author proves in \cite{\Efftpft}
that the synchronous and asynchronous versions of the \fftp\ are equivalent.

\section{The Example}
\begin{figure}[ht!]
  \begin{center}
  \psfrag{t}{\small$t$}
  \psfraga <-1pt,2pt> {a}{\small$a$}
  \includegraphics[width=9cm]{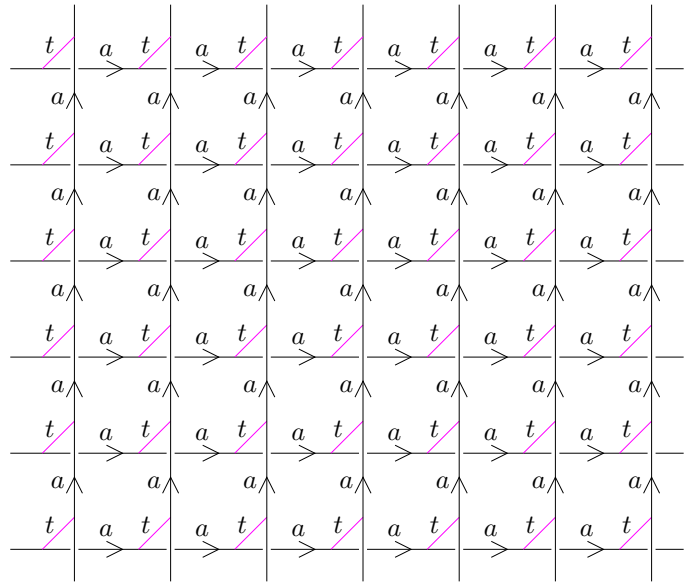}
  \end{center}
  \caption{The Cayley graph for $G=\langle a, t \; | \;  t^2=1,atat=tata
  \rangle$}
  \label{fig:cg}
\end{figure}

Let $G$ be the split extension of $\mathbb Z^2$, generated by
$\{a,b\}$, by $\mathbb Z_2$, generated by $\{t\}$, such that $t$
conjugates $a$ to $b$ and $b$ to $a$, with presentation
$$\langle a,b ,t \; | \;  t^2=1,ab=ba, tat=b\rangle.$$ Performing one Tietze
transformation (removing $b=tat$) we obtain
$$\langle a, t \; | \;  t^2=1,atat=tata \rangle.$$ Let $A=\{a^{\pm 1},t^{\pm
1}\}$ be the inverse-closed generating set corresponding to this
presentation.  The Cayley graph for $G$ with respect to $A$ is shown
in Figure \ref{fig:cg}.
We can consider the vertices of the graph as either being in the top
or the bottom layer, where edges labeled $t$ link top and bottom layers.
 We declare the identity vertex to lie in the
bottom. Each vertex can also be given a coordinate $(x,y)$ where $x$
is the distance in the East-West direction from the identity, and $y$
is the distance in the North-South direction. In this way each vertex
(group element) is uniquely specified by the triple
$(x,y,\mathrm{bottom})$ or $(x,y,\mathrm{top})$.

For example, the identity has the coordinate $(0,0,\mathrm{bottom})$,
the word $a^3ta^4$ has the coordinate $(3,4,\mathrm{top})$, the
word $a^3ta^4t$ has the coordinate $(3,4,\mathrm{bottom})$, and the
word $ta^3ta^4$ has the coordinate $(4,3,\mathrm{bottom})$.

\begin{lem}\label{lem:top=odd}
Each word from $1$ to a vertex in the top layer has an odd number of
$t$ letters, and each word from $1$ to a vertex in the bottom layer
has an even number of $t$ letters.
\end{lem}

\begin{proof}
Suppose a vertex has the coordinate $(x,y,\mathrm{top})$. Then
$a^xta^y$ is a word to this vertex. Suppose $w$ is any other word to
this vertex. Then $wa^{-y}ta^{-x}=_G 1$ so under the map which sends
$t$ to $t$ and $a$ to $1$ this word must be sent to an even power of
$t$, so $w$ has an odd number of $t$ letters. Similarly if a vertex
lies in the bottom layer there is a word $a^xta^yt$ to it from $1$. If
$w$ is any other word to this vertex then $wta^{-y}ta^{-x}=_G 1$ gets
sent to an even power of $t$ so $w$ has an even number of $t$ letters.
\end{proof}

\begin{lem}\label{lem:not-geod}
The word $ta^nta^mt$ for any $m,n\in \mathbb Z$ is not geodesic.
\end{lem}

\begin{proof}
The word $ta^nta^mt$ can be written as $a^mta^n$ which is shorter.
\end{proof}

\begin{lem}\label{lem:top=unique}
Each vertex in the top layer has a unique geodesic to it from the
identity of the form $a^xta^y$, where $(x,y,\mathrm{top})$ is the
coordinate of the vertex.
\end{lem}

\begin{proof}
By Lemma \ref{lem:top=odd} any geodesic to the vertex with coordinate
$(x,y,\mathrm{top})$ has an odd number of $t$ letters. If a word has
three or more $t$ letters then it has a subword of the form
$ta^nta^mt$ which is not geodesic by Lemma \ref{lem:not-geod}. So any
geodesic to this vertex has exactly one $t$ letter, so is of the form
$a^ita^j$. This path has coordinate $(i,j,\mathrm{top})$ so it must be
that $i=x$ and $j=y$.
\end{proof}

\begin{prop}\label{propn:notfftp}
$A$ does not have the \fftp.
\end{prop}

\begin{proof}
Suppose by way of contradiction that $A$ has the \fftp\ with positive
constant $k$, and choose $n >> k$.  Consider the word $w=ta^nta^nt$
which ends at the coordinate $(n,n,\mathrm{top})$. Any word that ends
at this vertex must move East at least $n$ units and North at least
$n$ units, and by Lemma \ref{lem:top=odd} must have an odd number of
$t$ letters. If it has just one $t$ letter then it must be the unique
geodesic $a^nta^n$ which clearly does not $k$-fellow travel
$w$. Otherwise it has three or more $t$ letters, so has length at
least $2n+3$ so is not shorter than $w$, so we are done.
\end{proof}

\begin{thm}\label{thm:main}
There is a group and finite \gset\ such that the language of all
geodesics is regular but fails to have the \fftp.
\end{thm}
\begin{proof}
By Proposition \ref{propn:notfftp} the group $G$ with \gset\ $A$ fails
the \fftp.  

Consider the language $L=\{a^x,a^xta^y,a^{x_1}ta^yta^{x_2}:x,x_1,x_2,y
\in \mathbb Z, x_1\cdot x_2\geq0\}$.  $L$ is the language of the
finite state automaton in Figure \ref{fig:fsa}. All states are accept
states.
\begin{figure}[ht!]
  \begin{center}
  \psfrag{a}{\small$a$}
  \psfraga <-3pt,-2pt> {ai}{\small$a^{{-}1}$}
  \psfraga <-3pt,-8pt> {aii}{\small$a^{{-}1}$}
  \psfraga <-4pt,1pt> {s}{\small$s_0$}
  \includegraphics[width=11cm]{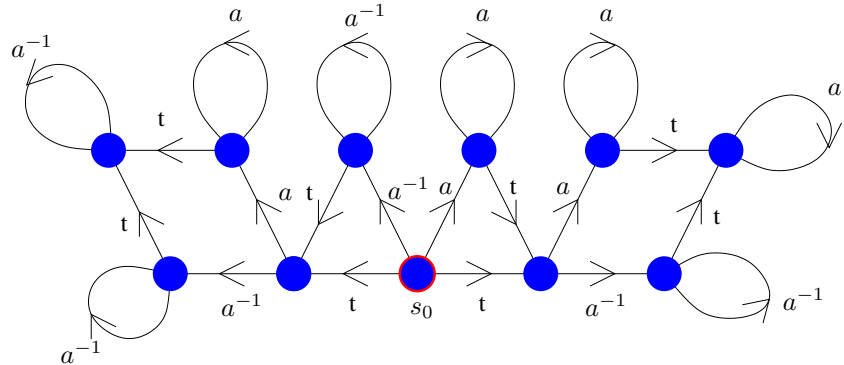}
  \end{center}
  \caption{A finite state automaton accepting the language $L$}
  \label{fig:fsa}
\end{figure}

We now show that $L$ is the language of all geodesics on $A$ for $G$.
By Lemma \ref{lem:top=unique} every group element corresponding to a
vertex in the top layer has a unique geodesic representative of the
form $a^xta^y$. Otherwise the group element corresponds to a vertex in
the bottom, so has an even number of $t$ letters.  If a word has more
than two $t$ letters then it has a subword of the form $ta^nta^mt$
which is not geodesic by Lemma \ref{lem:not-geod}, so a geodesic word
for a bottom element has either zero $t$ letters, so is of the form
$a^x$, or has two $t$ letters, so is of the form
$a^{x_1}ta^yta^{x_2}$. If $x_1$ and $x_2$ don't have the same sign
then we can find a shorter word $a^{(x_1+x_2)}ta^yt$. Otherwise
$a^{x_1}ta^yta^{x_2}$ is a geodesic to a vertex with coordinate
$(x_1+x_2,y,\mathrm{bottom})$. Notice that this gives a family of
$x_1+x_2+1$ geodesics to this vertex.
\end{proof}

\Addresses\recd

\end{document}